\documentclass[10pt]{amsart}

\usepackage[latin1]{inputenc}
\usepackage{amsmath}
\usepackage{amsfonts}
\usepackage{amssymb}
\usepackage[mathscr]{eucal}
\usepackage{diagrams}
\usepackage{xy,color,graphicx}
\usepackage{hyperref}
\xyoption{all}

\newcommand{\wis}[1]{{\text{\em \usefont{OT1}{cmtt}{m}{n} #1}}}

\newcommand{\PP}{\mathbb{P}}
\newcommand{\C}{\mathbb{C}}
\newcommand{\Z}{\mathbb{Z}}
\newcommand{\vtx}[1]{*+[o][F-]{\scriptscriptstyle #1}}

\newtheorem{proposition}{Proposition}

\newtheorem{theorem}{Theorem}
\newtheorem{corollary}{Corollary}
\newtheorem{lemma}{Lemma}

\newtheorem{remark}{Remark}

\title{The geometry of representations of \\ 3-dimensional Sklyanin algebras}
\author{Kevin De Laet}
\address{Mathematics and statistics, Hasselt University \\
Agoralaan - Building D, B-3590 Diepenbeek (Belgium) \\ {\tt kevin.delaet@uhasselt.be}}
\author{Lieven Le Bruyn} 
\address{Department of Mathematics, University of Antwerp \\ 
 Middelheimlaan 1, B-2020 Antwerp (Belgium) \\ {\tt lieven.lebruyn@uantwerpen.be}}

\begin{document}
\sloppy

\maketitle

\begin{abstract}
The representation scheme $\wis{rep}_n~A$ of the 3-dimensional Sklyanin algebra $A$ associated to a plane elliptic curve and $n$-torsion point contains singularities over the augmentation ideal $\mathfrak{m}$. We investigate the semi-stable representations of the noncommutative blow-up algebra $B=A \oplus \mathfrak{m}t \oplus \mathfrak{m}^2 t^2 \oplus \hdots$ to obtain a partial resolution of the central singularity
\[
\wis{proj}~Z(B) \rOnto \wis{spec}~Z(A) \]
such that the remaining singularities in the exceptional fiber determine an elliptic curve and are all of type $\C \times \C^2/\Z_n$.
\end{abstract}

\section{Introduction}

Three dimensional Sklyanin algebras appear in the classification by M. Artin and W. Schelter \cite{AS} of graded algebras of global dimension $3$. In the early 90ties this class of algebras was studied extensively by means of noncommutative projective algebraic geometry, see a.o. \cite{Artin}, \cite{ATV1}, \cite{ATV2}, \cite{LBsymbols} and \cite{TateSmith}. Renewed interest in this class of algebras arose recently as they are superpotential algebras and as such relevant in supersymmetric quantum field theories, see a.o. \cite{Berenstein} and \cite{Walton}.

Consider a smooth elliptic curve $E$ in Hesse normal form $\mathbb{V}((a^3+b^3+c^3)XYZ-abc(X^3+Y^3+Z^3)) \rInto \mathbb{P}^2$ and the point $p=[a:b:c]$ on $E$. The $3$-dimensional Sklyanin algebra $A$ corresponding to the pair $(E,p)$ is the noncommutative algebra with defining equations
\[
\begin{cases}
a xy + b yx + c z^2 = 0 \\
a yz + b zy + c x^2 = 0 \\
a zx + b xz + c y^2 = 0
\end{cases}
\]
The connection comes from the fact that the multi-linearization of these equations defines a closed subscheme in $\mathbb{P}^2 \times \mathbb{P}^2$ which is the graph of translation by $p$ on the elliptic curve $E$, see \cite{ATV1}. Alternatively, one obtains the defining equations of $A$ from the superpotential $W=a xyz + b yzx + \frac{c}{3}(x^3+y^3+z^3)$, see \cite{Walton}.

The algebra $A$ has a central element of degree $3$, found by computer search in \cite{AS}
\[
c_3 = c(a^3-c^3)x^3+a(b^3-c^3)xyz + b(c^3-a^3)yxz + c(c^3-b^3)y^3 \]
with the property that $A/(c_3)$ is the twisted coordinate ring of the elliptic curve $E$ with respect to the automorphism given by translation by $p$, see \cite{ATV1}. We will prove an intrinsic description of this central element, answering a MathOverflow question \cite{euclid}.

\begin{theorem}
The central element $c_3$ of the $3$-dimensional Sklyanin algebra $A$ corresponding to the pair $(E,p)$ can be written as
\[
c(a^3-b^3)(xyz+yzx+zxy)+b(c^3-a^3)(yxz+xzy+zyx)+c(a^3-b^3)(x^3+y^3+zx^3) \]
and is the superpotential of the $3$-dimensional Sklyanin algebra $A'$ corresponding to the pair $(E,[-2]p)$.
\end{theorem}

Next, we turn to the study of finite dimensional representations of $A$ which is important in supersymmetric gauge theory as they correspond to the vacua states. It is well known that $A$ is a finite module over its center $Z(A)$ and a maximal order in a central simple algebra of dimension $n^2$ if and only if the point $p$ is of finite order $n$, see \cite{ATV1}. We will further assume that $(n,3)=1$ in which case J. Tate and P. Smith proved in \cite{TateSmith} that the center $Z(A)$ is generated by $c_3$ and the reduced norms of $x,y$ and $z$ (which are three degree $n$ elements, say $x',y',z'$) satisfying one relation of the form
\[
c_3^n = cubic(x',y',z') \]
It is also known that $\wis{proj}~Z(A) \simeq \mathbb{P}^2$ with coordinates $[x':y':z']$ in which the $cubic(x',y',z')$ defines the isogenous elliptic curve $E'=E/\langle p \rangle$, see a.o. \cite{LBsymbols}. We will use these facts to give explicit matrices for the simple $n$-dimensional representations of $A$ and show that $A$ is an Azumaya algebra away from the isolated central singularity.

However, the scheme $\wis{rep}_n~A$ of all (trace preserving) $n$-dimensional representations of $A$ contains singularities in the nullcone. We then try to resolve these representation singularities by considering the noncommutative analogue of a blow-up algebra
\[
B = A \oplus \mathfrak{m} t \oplus \mathfrak{m}^2 t^2 \oplus \hdots \subset A[t,t^{-1}] \]
where $\mathfrak{m}=(x,y,z)$ is the augmentation ideal of $A$. We will prove

\begin{theorem} The scheme $\wis{rep}^{ss}_n~B$ of all semi-stable $n$-dimensional representations of the blow-up algebra $B$ is a smooth variety.
\end{theorem}

This allows us to compute all the (graded) local quivers in the closed orbits of $\wis{rep}^{ss}_n~B$ as in \cite{LBBook} and \cite{BocklandtSymens}. This information then leads to the main result of this paper which gives a partial resolution of the central isolated singularity.

\begin{theorem} The exceptional fiber $\mathbb{P}^2$ of the canonical map
\[
\wis{proj}~Z(B) \rOnto \wis{spec}~Z(A) \]
contains $E' = E/ \langle p \rangle$ as the singular locus of $\wis{proj}~Z(B)$. Moreover, all these singularities are of type $\C \times \C^2/\mathbb{Z}_n$ with $\C^2/\mathbb{Z}_n$ an Abelian quotient surface singularity.
\end{theorem}

\section{Central elements and superpotentials}

The finite Heisenberg group of order $27$
\[
\langle~u,v,w~|~[u,v]=w,~[u,w]=[v,w]=1,~u^3=v^3=w^3=1~\rangle \]
has a $3$-dimensional irreducible representation $V = \C x + \C y + \C z$ given by the action
\[
u \mapsto \begin{bmatrix} 0 & 1 & 0 \\ 0 & 0 & 1 \\ 1 & 0 & 0 \end{bmatrix} \quad
v \mapsto \begin{bmatrix} 1 & 0 & 0 \\ 0 & \rho & 0 \\ 0 & 0 & \rho^2 \end{bmatrix} \quad
w \mapsto \begin{bmatrix} \rho & 0 & 0 \\ 0 & \rho & 0 \\ 0 & 0 & \rho \end{bmatrix} \]
One verifies that $V \otimes V$ decomposes as three copies of $V^*$, that is,
\[
V \otimes V \simeq \wedge^2(V) \oplus S^2(V) \simeq V^* \oplus (V^* \oplus V^*) \]
where the three copies can be taken to be the subspaces
\[
\begin{cases}
V_1 = \C (yz-zy) + \C (zx-xz) + \C (xy-yx) \\
V_2 = \C (yz+zy) + \C (zx+xz) + \C (xy+yx) \\
V_3 = \C x^2 + \C y^2 + \C z^2
\end{cases}
\]
Taking the quotient of $\C \langle x,y,z \rangle$ modulo the ideal generated by $V_1 = \wedge^2 V$ gives the commutative polynomial ring $\C[x,y,z]$. Hence we can find analogues of the polynomial ring in three variables by dividing $\C \langle x,y,z \rangle$ modulo the ideal generated by another copy of $V^*$ in $V \otimes V$ and the resulting algebra will inherit an action by $H_3$. Such a copy of $V^*$ exists for all $[A:B:C] \in \PP^2$ and is spanned by the three vectors
\[
\begin{cases}
A(yz-zy) + B (yz+zy) + C x^2 \\
A(zx-xz) + B (zx+xz) + C y^2 \\
A(xy-yx) + B (xy+yx) + C z^2
\end{cases}
\]
and by taking $a=A+B, b=B-A$ and $c=C$ we obtain the defining relations of the $3$-dimensional Sklyanin algebra. In particular there is an $H_3$-action on $A$ and the canonical central element $c_3$ of degree $3$ must be a 1-dimensional representation of $H_3$. It is obvious that $c_3$ is fixed by the action of $v$ and a minor calculation shows that $c_3$ is also fixed by $u$. Therefore, the central element $c_3$ given above, or rather $3c_3$, can also be represented as
\[
a(b^3-c^3)(xyz+yzx+zxy)+b(c^3-a^3)(yxz+xzy+zyx)+c(a^3-b^3)(x^3+y^3+z^3) \]
Now, let us reconsider the superpotential $W=a xyz + b yxz + \frac{c}{3}(x^3+y^3+z^3)$ for a $[a:b:c] \in \mathbb{P}^2$. This superpotential gives us three quadratic relations by taking cyclic derivatives with respect to the variables
\[
\begin{cases}
\partial_x W = a yz + b z y + c x^2 \\
\partial_y W = a zx + b xz + c y^2 \\
\partial_z W = a xy + b yx + c z^2
\end{cases}
\]
giving us the defining relations of the $3$-dimensional Sklyanin algebra. We obtain the same equations by considering a more symmetric form of $W$, or rather of $3W$
\[
a (xyz+yzx+zxy) + b (yxz+xzy+zyx) + c (x^3+y^3+z^3) \]
We see that the form of the central degree $3$ element and of the superpotential are similar but with different coefficients. This means that the central element is the superpotential defining another $3$-dimensional Sklyanin algebra and theorem~1 clarifies this connection.

\vskip 3mm
\noindent
{\bf Proof of Theorem~1 :} The $3$-dimensional Sklyanin algebra determined by the superpotential $3c_3$ is determined by $[a(b^3-c^3) :b(c^3-a^3):c(a^3-b^3)]$ (instead of $[a:b:c]$ for the original). Therefore, the associated elliptic curve has defining Hesse equation
\[
\mathbb{V}(\alpha (x^3+y^3+z^3) - \beta xyz) \rInto \mathbb{P}^2 \]
where
\[
\begin{cases}
\alpha = a(b^3-c^3)b(c^3-a^3)c(a^3-b^3) \\
\beta = (a(b^3-c^3))^3+(b(c^3-a^3))^3+(c(a^3-b^3))^3
\end{cases}
\]
but this is the same equation, upto a scalar, as the original curve
\[
E = \mathbb{V}(abc(x^3+y^3+c^3) - (a^3+b^3+c^3)xyz) \]
The tangent line to $E$ in the point $p=[a:b:c]$ has equation
\[
\mathbb{V}((2a^3bc-b^4c-bc^4)(x-a)+(2ab^3c-a^4c-ac^4)(y-b)+(2abc^3-a^4b-ab^4)(z-c)) \]
and so the third point of intersection is
\[
[-2]p = [a(b^3-c^3):b(c^3-a^3):c(a^3-b^3)] \]
which are the parameters of the algebra.  $\hfill \qed$

\section{Resolving representation singularities}

Let $R$ be a graded $\C$-algebra, generated by finitely many elements $x_1,\hdots,x_m$ where $deg(x_i)=d_i \geq 0$, which is a finite module over its center $Z(R)$. Following \cite{ProcesiCH} we say that $R$ is a Cayley-Hamilton algebra of degree $n$ if there is a $Z(R)$-linear gradation preserving trace map
$tr~:~R \rTo Z(R)$ such that for all $a,b \in R$ we have
\begin{itemize}
\item{$tr(ab)=tr(ba)$}
\item{$tr(1)=n$}
\item{$\chi_{n,a}(a) = 0$}
\end{itemize}
where $\chi_{n,a}(t)$ is the $n$-th Cayley-Hamilton identity expressed in the traces of powers of $a$. Maximal orders in a central simple algebra of dimension $n^2$ are examples of Cayley-Hamilton algebras of degree $n$. 

In particular, a $3$-dimensional Sklyanin algebra $A$ associated to a couple $(E,p)$ where $p$ is a torsion point of order $n$, and the corresponding blow-up algebra $B = A \oplus \mathfrak{m} t \oplus \mathfrak{m}^2 t^2 \oplus \hdots$ are affine graded Cayley-Hamilton algebras of degree $n$ equipped with the (gradation preserving) reduced trace map.

If $R$ is an affine graded Cayley-Hamilton algebra of degree $n$ we define $\wis{rep}_n~R$ to be the affine scheme of all $n$-dimensional trace preserving representations, that is of all algebra morphisms
\[
R \rTo^{\phi} M_n(\C) \qquad \text{such that} \qquad \forall a \in R~:~\phi(tr(a)) = Tr(\phi(a)) \]
where $Tr$ is the usual trace map on $M_n(\C)$. Isomorphism of representations defines a $\wis{GL}_n$-action of $\wis{rep}_n~R$ and a result of Artin's \cite{Artin69} asserts that the closed orbits under this action, that is the points of the GIT-quotient scheme $\wis{rep}_n~R//\wis{GL}_n$,  are precisely the isomorphism classes of $n$-dimensional trace preserving semi-simple representations of $R$. The reconstruction result of Procesi \cite{ProcesiCH} asserts that in this setting
\[
\wis{spec}~Z(R) \simeq \wis{rep}_n~R // \wis{GL}_n \]
The gradation on $R$ defines an additional $\C^*$-action on $\wis{rep}_n~R$ commuting with the $\wis{GL}_n$-action. With $\wis{rep}_n^{ss}~R$ we denote the Zariski open subset of all semi-stable trace preserving representations $\phi~:~R \rTo M_n(\C)$, that is, such that there is an homogeneous central element $c$ of positive degree such that $c(\phi) \not= 0$. We have the following graded version of Procesi's reconstruction result, see a.o. \cite{BocklandtSymens}
\[
\wis{proj}~Z(R) \simeq \wis{rep}^{ss}_n~R // \wis{GL}_n \times \C^* \]
As a $\wis{GL}_n \times \C^*$-orbit is closed in $\wis{rep}_n^{ss}~R$ if and only if the $\wis{GL}_n$-orbit is closed we see that points of $\wis{proj}~Z(R)$ classify one-parameter families of isoclasses of trace-preserving $n$-dimensional semi-simple representations of $R$. In case of a simple representation such a one-parameter family determines a graded algebra morphism
\[
R \rTo M_n(\C[t,t^{-1}])(\underbrace{0,\hdots,0}_{m_0},\underbrace{1,\hdots,1}_{m_1},\hdots,\underbrace{e-1,\hdots,e-1}_{m_{e-1}}) \]
where $e$ is the degree of $t$ and where we follow \cite{NastaFred} in defining the shifted graded matrix algebra $M_n(S)(a_1,\hdots,a_n)$ by taking is homogeneous part of degree $i$ to be
\[
\begin{bmatrix}
S_i & S_{i-a_1+a_2} & \hdots & S_{i-a_1+a_n} \\
S_{i-a_2+a_1} & S_i & \hdots & S_{i-a_2+a_n} \\
\vdots & \vdots & \ddots & \vdots \\
S_{i-a_n+a_1} & S_{i-a_n+a_2} & \hdots & S_i \end{bmatrix} \]
The $\wis{GL}_n \times \C^*$-stabilizer subgroup of any of the simples $\phi$ in this family is then isomorphic to $\C^* \times \pmb{\mu}_e$ where the cyclic group $\pmb{\mu}_e$ has generator $(g_{\zeta},\zeta) \in \wis{GL}_n \times \C^*$ where $\zeta$ is a primitive $e$-th root of unity and
\[
g_{\zeta} = \wis{diag}(\underbrace{1,\hdots,1}_{m_0},\underbrace{\zeta,\hdots,\zeta}_{m_1},\hdots,\underbrace{\zeta^{e-1},\hdots,\zeta^{e-1}}_{m_{e-1}}) \]
see \cite[lemma 4]{BocklandtSymens}. If, in addition, $\phi$ is a smooth point of $\wis{rep}^{ss}_n~R$ then the normal space 
\[
N(\phi) = T_{\phi} \wis{rep}_n^{ss}R / T_{\phi} \wis{GL}_n.\phi \]
to the $\wis{GL}_n$-orbit decomposes as a $\pmb{\mu}_e$-representation into a direct sum of $1$-dimensional simples
\[
N(\phi) = \C_{{i_1}} \oplus \hdots \oplus \C_{{i_d}} \]
where the action of the generator on $\C_{_k}$ is by multiplication with $\zeta^k$. Alternatively, $\phi$ determines a (necessarily smooth) point $[\phi] \in \wis{spec}~Z(R)$ and because $N(\phi)$ is equal to $Ext^1_R(S_{\phi},S_{\phi})$ and because $R$ is Azumaya in $[\phi]$ it coincides with $Ext^1_{Z(R)}(S_{[\phi]},S_{[\phi]})$ (where $S_{[\phi]}$ is the simple $1$-dimensional representation of $Z(R)$ determined by $[\phi]$) which is identical to the tangent space $T_{[\phi]} \wis{spec} Z(R)$. The action of the stabilizer subgroup $\pmb{\mu}_e$ on $Ext^1_R(S_{\phi},S_{\phi})$ carries over to that on $T_{[\phi]} \wis{spec} Z(R)$. 

The one-parameter family of simple representations also determines a point $\overline{\phi} \in \wis{proj} Z(R)$ and an application of the Luna slice theorem \cite{Luna} asserts that for all $t \in \C$ there is a neighborhood of $(\overline{\phi},t) \in \wis{proj} Z(R) \times \C$ which is \'etale isomorphic to a neighborhood of $0$ in $N(\phi) // \pmb{\mu}_e$, see \cite[Thm. 5]{BocklandtSymens}. 

\subsection{From $\wis{Proj}(A)$ to $\wis{rep}_n A$} In noncommutative projective algebraic geometry, see a.o. \cite{ATV1},\cite{ATV2} and \cite{Artin}, one studies the Grothendieck category $\wis{Proj}(A)$ which is the quotient category of all graded left $A$-modules modulo the subcategory of torsion modules. In the case of $3$-dimensional Sklyanin algebras the linear modules, that is those with Hilbert series $(1-t)^{-1}$ (point modules) or $(1-t)^{-2}$ (line modules) were classified in \cite{ATV2}. Identify $\mathbb{P}^2$ with $\mathbb{P}^2_{nc} = \mathbb{P}(A_1^*)$, then
\begin{itemize}
\item{point modules correspond to points on the elliptic curve $E \rInto \mathbb{P}^2_{nc}$}
\item{line modules correspond to lines in $\mathbb{P}^2_{nc}$}
\end{itemize}
In the case of interest to us, when $A$ corresponds to a couple $(E,p)$ with $p$ a torsion point of order $n$ also fat modules are important which are critical cyclic graded left $A$-modules with Hilbert series $n.(1-t)^{-1}$. They were classified by M. Artin \cite{Artin} and are relevant in the study of $\wis{proj} Z(A) = \mathbb{P}^2_c = \mathbb{P}(Z(A)_n^*)$. Observe that the reduced norm map $N$ relates the different manifestations of $\mathbb{P}^2$ and the elliptic curve $E$ with its isogenous curve $E/ \langle p \rangle$
\[
\xymatrix{\PP^2_{nc}=\PP(A_1^*) \ar[rr]^{N} & & \PP^2_c = \PP(Z(A)_n^*) \\
E \ar[u] \ar[rr]^{./ \langle p \rangle} & & E' = E/ \langle p \rangle \ar[u]}
\]
Points $\rho \in \mathbb{P}^2_c - E'$ determine fat points $F_{\pi}$ with graded endomorphism ring isomorphic to $M_n(\C[t,t^{-1}])$ with $deg(t)=1$, and hence determine a one-parameter family of simple $n$-dimensional representations in $\wis{rep}^{ss}_n A$ with $\wis{GL}_n \times \C^*$-stabilizer subgroup $\C^* \times 1$. There is an effective method to construct $F_{\pi}$, see \cite{LBsymbols}. Write $\rho$ as the intersection of two lines $\mathbb{V}(z) \cap \mathbb{V}(z')$ and let $\mathbb{V}(z') \cap E' = \{ q_1,q_2,q_3 \}$ be the intersection with the elliptic curve $E'$. Then by lifting the $q_i$ through the isogeny to $n$ points $p_{ij} \in E$ we see that we can lift the line $\mathbb{V}(z')$ to $n^2$ lines in $\mathbb{P}^2_{nc} = \mathbb{P}(A_1^*)$, that is, there are $n^2$ one-dimensional subspaces $\C l \subset A_1$ with the property that $\C N(l) = \C z'$. The fat point corresponding to $\pi$ is then the shifted quotient of a line module determined by $l$
\[
F_{\rho} \simeq \frac{A}{A.l+A.z}[n] \]
On the other hand, if $q$ is a point on $E'$, then lifting $q$ through the isogeny results in an orbit of $n$ points of $E$, $\{ r, r+p, r+[2]p, \hdots, r+[n-1]p \}$. If $P$ is the point module corresponding to $r \in E$, then the fat point module corresponding to $q$ is
\[
F_q = P \oplus P[1] \oplus P[2] \oplus \hdots \oplus P[n-1] \]
and the corresponding graded endomorphism ring is isomorphic to $M_n(\C[t,t^{-1}])(0,1,2,\hdots,n-1)$ where $deg(t)=n$ and hence corresponds to a one-parameter family of simple $n$-dimensional representations in $\wis{rep}^{ss}_n A$ with $\wis{GL}_n \times \C^*$-stabilizer subgroup generated by $\C^* \times 1$ and a cyclic group of order $n$
\[
\pmb{\mu}_n = \langle ( \begin{bmatrix} 1 & & & \\ & \zeta & & \\ & & \ddots & \\ & & & \zeta^{n-1} \end{bmatrix}, \zeta ) \rangle \]
with $\zeta$ a primitive $n$-th root of unity. In fact, we can give a concrete matrix-representation of these simple modules. Assume that $r+[i]p = [a_i:b_i:c_i] \in \mathbb{P}^2_{nc}$ then the fat point module $F_q$ corresponds to the quiver-representation
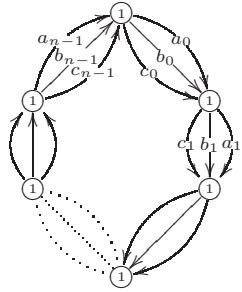
\[
\xymatrix{
& \vtx{1} \ar@/^2ex/[rd]|{a_0} \ar[rd]|{b_0} \ar@/_2ex/[rd]|{c_0} & \\
\vtx{1} \ar@/^2ex/[ru]|{a_{n-1}} \ar[ru]|{b_{n-1}} \ar@/_2ex/[ru]|{c_{n-1}}& & \vtx{1} \ar@/^2ex/[d]|{a_1} \ar[d]|{b_1} \ar@/_2ex/[d]|{c_1}\\
\vtx{1} \ar@/^2ex/[u] \ar[u] \ar@/_2ex/[u] & & \vtx{1} \ar@/^2ex/[ld] \ar[ld] \ar@/_2ex/[ld] \\
& \vtx{1} \ar@{.}@/^2ex/[lu] \ar@{.}[lu] \ar@{.}@/_2ex/[lu] & }
\]
and the map $A \rTo M_n(\C[t,t^{-1}])(0,1,2,\hdots,n-1)$ sends the generators $x,y$ and $z$ to the degree one matrices
{\tiny
\[
\begin{bmatrix}
0 & 0 &  \hdots & \hdots & a_{n-1} t \\
a_0 & 0 &  \hdots & \hdots & 0 \\
0 & a_1 &  \ddots &  & \vdots \\
\vdots &  & \ddots & \ddots & \vdots \\
0 & 0  & \hdots & a_{n-2} & 0 
\end{bmatrix} \quad 
\begin{bmatrix}
0 & 0 &  \hdots & \hdots & b_{n-1} t \\
b_0 & 0 &  \hdots & \hdots & 0 \\
0 & b_1 &  \ddots &  & \vdots \\
\vdots &  & \ddots & \ddots & \vdots \\
0 & 0  & \hdots & b_{n-2} & 0 
\end{bmatrix}
\quad 
\begin{bmatrix}
0 & 0 &  \hdots & \hdots & c_{n-1} t \\
c_0 & 0 &  \hdots & \hdots & 0 \\
0 & c_1 &  \ddots &  & \vdots \\
\vdots &  & \ddots & \ddots & \vdots \\
0 & 0  & \hdots & c_{n-2} & 0 
\end{bmatrix}
\]
}

\begin{theorem} Let $A$ be a $3$-dimensional Sklyanin algebra corresponding to a couple $(E,p)$ where $p$ is a torsion point of order $n$ and assume that $(n,3)=1$. Consider the GIT-quotient
\[
\wis{rep}_n A \rOnto^{\pi} \wis{spec} Z(A) = \wis{rep}_n A // \wis{GL}_n \]
Then we have
\begin{enumerate}
\item{$\wis{rep}_n^{ss} A$ is a smooth variety of dimension $n^2+2$}
\item{$A$ is an Azumaya algebra away from the isolated singularity $\tau \in \wis{spec} Z(A)$}
\item{the nullcone $\pi^{-1}(\tau)$ contains singularities}
\end{enumerate}
\end{theorem}

\begin{proof} We know that $\mathbb{P}^2_c = \wis{proj} Z(A) = \wis{rep}^{ss}_n A // \wis{GL}_n \times \C^*$ classifies one-parameter families of semi-stable $n$-dimensional semi-simple representations of $A$. To every point $\rho \in \mathbb{P}^2_c$ we have associated a one-parameter family of simples, so all semi-stable $A$-representations are in fact simple as the semi-simplification $M^{ss}$ of a semi-stable representation still belongs to $\wis{rep}^{ss}_n A$. But then,
all non-trivial semi-simple $A$-representations are simple and therefore the GIT-quotient
\[
\wis{rep}_n^{ss} A \rOnto \wis{spec} Z(A) - \{ \tau \} = \wis{rep}^{ss}_n A// \wis{GL}_n \]
is a principal $\wis{PGL}_n$-fibration in the \'etale topology. This proves (1).

The second assertion follows as principal $\wis{PGL}_n$-fibrations in the \'etale topology correspond to Azumaya algebras. For (3), if $\wis{rep}_n A$ would be smooth, the algebra $A$ would be Cayley-smooth as in \cite{LBBook}. There it is shown that the only type of central singularity that can arise for Cayley-smooth algebras with a $3$-dimensional center is the conifold singularity.
\end{proof}

If we want to distinguish between the two types of simple representations, we have to consider the $\wis{GL}_n \times \C^*$-action.

\begin{lemma} If $S$ is a simple $A$-representation with $\wis{GL}_n \times \C^*$-orbit determining a fat point $F_q$ with $q \in E'$, then the normal space $N(S)$ to the $\wis{GL}_n$-orbit decomposes as representation over the $\wis{GL}_n \times \C^*$-stabilizer subgroup $\pmb{\mu}_n$ as $\C_0 \oplus \C_0 \oplus \C_3$, or in the terminology of \cite{BocklandtSymens}, the associated local weighted quiver is
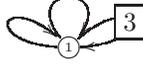
\[
\xymatrix{ \vtx{1} \ar@(l,ul) \ar@(ul,ur) \ar@(ur,r)|{\boxed 3}} \]
\end{lemma}

\begin{proof}
From \cite{TateSmith} we know that the center $Z(A)$ can be represented as
\[
Z(A) = \frac{\C[x',y',z',c_3]}{(c_3^n - cubic(x',y',z'))} \]
where $x',y',z'$ are of degree $n$ (the reduced norms of $x,y,z$) and $c_3$ is the canonical central element of degree $3$. The simple $A$-representation $S$ determines a point $s \in \wis{spec} Z(A)$ such that $c_3(s)=0$. Again, as $A$ is Azumaya over $s$ we have that $N(S)=Ext^1_A(S,S)$ coincides with the tangent space $T_s \wis{spec} Z(A)$. Gradation defines a $\pmb{\mu}_n$-action on $Z(A)$ leaving $x',y',z'$ invariant and sending $c_3$ to $\zeta^3 c_3$. The stabilizer subgroup of this action in $s$ is clearly $\pmb{\mu}_n$ and computing the tangent space gives the required decomposition.
\end{proof}

\subsection{$\mathcal{A}$ is Cayley-smooth}
Because $A$ is a finitely generated module over $Z(A)$, it defines a coherent sheaf of algebras $\mathcal{A}$ over $\wis{proj} Z(A) = \PP^2$. In this subsection we will show that $\mathcal{A}$ is a sheaf of Cayley-smooth algebras of degree $n$.

As $(n,3)=1$ it follows that the graded localisation $Q^g_{x'}(A)$ at the multiplicative set of central elements $\{ 1,x',x'^2,\hdots \}$ contains central elements $t$ of degree one and hence is isomorphic as a graded algebra to
\[
Q^g_{x'}(A) = (Q^g_{x'}(A))_0[t,t^{-1}] \]
By definition $\Gamma(\mathbb{X}(x'),\mathcal{A}) = (Q^g_{x'}(A))_0$ and by the above isomorphism it follows that $\Gamma(\mathbb{X}(x'),\mathcal{A})$ is a Cayley-Hamilton domain of degree $n$ and is Auslander regular of dimension two and consequently a maximal order. Repeating this argument for the other standard opens $\mathbb{X}(y')$ and $\mathbb{X}(z')$ we deduce

\begin{proposition} 
$\mathcal{A}$ is a coherent sheaf of Cayley-Hamilton maximal orders of degree $n$ which are Auslander regular domains of dimension $2$ over $\wis{proj} Z(A) = \PP^2_c$.
\end{proposition}

Thus, $\mathcal{A}$ is a maximal order over $\PP^2$ in a division algebra $\Sigma$ over $\C(\PP^2)$ of degree $n$. By the Artin-Mumford exact sequence (see for example \cite[3.6]{LBBook}) describing the Brauer group of $\C(\PP^2)$ we know that $\Sigma$ is determined by the ramification locus of $\mathcal{A}$ together with a cyclic $\Z_n$-cover over it.

Again using the above local description of $A$ as a graded algebra over $Z(A)$ we see that the fat point module corresponding to a point $p \notin E'$ determines a simple $n$-dimensional representation of $\mathcal{A}$ and therefore $\mathcal{A}$ is Azumaya in $p$. However, if $p \in E'$, then the corresponding fat point is of the form $P \oplus P[1] \oplus \hdots \oplus P[n-1]$ and this corresponds to a semi-simple $n$-dimensional representation which is the direct sum of $n$ distinct one-dimensional $\mathcal{A}$-representations, one component for each point of $E$ lying over $p$. Hence, we see that the ramification divisor of $\mathcal{A}$ coincides with $E'$ and, naturally, the division algebra $\Sigma$ is the one coresponding to the cyclic $\Z_n$-cover $E \rOnto E' = E/\langle \tau \rangle$.

Because $\mathcal{A}$ is a maximal order with smooth ramification locus, we deduce from \cite[\S 5.4]{LBBook}

\begin{proposition} $\mathcal{A}$ is a sheaf of Cayley-smooth algebras over $\PP^2_c$ and hence $\wis{rep}_n(\mathcal{A})$ is a smooth variety of dimension $n^2+1$ with GIT-quotient
\[
\wis{rep}_n(\mathcal{A}) \rOnto^{\pi} \PP^2_c = \wis{rep}_n(\mathcal{A})//\wis{GL}_n \]
and is a principal $\wis{PGL}_n$-fibration over $\PP^2_c - E'$. 
\end{proposition}

\subsection{The non-commutative blow-up} Consider the augmentation ideal $\mathfrak{m} = (x,y,z)$ of the $3$-dimensional Sklyanin algebra $A$ corresponding to a couple $(E,p)$ with $p$ a torsion point of order $n$. Define the non-commutative blow-up algebra to be the graded algebra
\[
B = A \oplus \mathfrak{m} t \oplus \mathfrak{m}^2 t^2 \oplus \hdots \subset A[t] \]
with degree zero part $A$ and where the commuting variable $t$ is given degree $1$. Note that $B$ is a graded subalgebra of $A[t]$ and therefore is again a Cayley-Hamilton algebra of degree $n$. Moreover, $B$ is a finite module over its center $Z(B)$ which is a graded subalgebra of $Z(A)[t]$.
Observe that $B$ is generated by the degree zero elements $x,y,z$ and by the degree one elements $X=xt,Y=yt$ and $Z=zt$. Apart from the Sklyanin relations among $x,y,z$ and among $X,Y,Z$ these generators also satisfy commutation relations such as $Xx=xX,Xy=xY,Xz=xZ$ and so on.

With $\wis{rep}^{ss}_n B$ we will denote again the Zariski open subset of $\wis{rep}_n B$ consisting of all trace-preserving $n$-dimensional semi-stable representations, that is, those on which some central homogeneous element of $Z(B)$ of strictly positive degree does not vanish. Theorem 2 asserts that $\wis{rep}^{ss}_n B$ is a smooth variety of dimension $n^2+3$.

\vskip 3mm
\noindent
{\bf Proof of Theorem 2 : } As before, we have a $\wis{GL}_n \times \C^*$-action on $\wis{rep}_n^{ss} B$ with corresponding GIT-quotient
\[
\wis{proj} Z(B) \simeq \wis{rep}^{ss}_n B // \wis{GL}_n \times \C^* \]
Composing the GIT-quotient map with the canonical morphism (taking the degree zero part) $\wis{proj} Z(B) \rOnto \wis{spec} Z(A)$ we have a projection
\[
\gamma~:~\wis{rep}^{ss}_n B \rOnto \wis{spec} Z(A) \]
Let $\mathfrak{p}$ be a maximal ideal of $Z(A)$ corresponding to a smooth point, then the graded localization of $B$ at the degree zero multiplicative subset $Z(A)-\mathfrak{p}$ gives
\[
B_{\mathfrak{p}} \simeq A_{\mathfrak{p}}[t,t^{-1}] \]
whence $B_{\mathfrak{p}}$ is an Azumaya algebra over $Z(A)[t,t^{-1}]$ and therefore over $\wis{spec} Z(A) - \{ \tau \}$ the projection $\gamma$ is a principal $\wis{PGL}_n \times \C^*$-fibration and in particular the dimension of $\wis{rep}^{ss}_n B$ is equal to $n^2+3$.

This further shows that possible singularities of $\wis{rep}^{ss}_n B$ must lie in $\gamma^{-1}(\tau)$ and as the singular locus is Zariski closed we only have to prove smoothness in points of closed $\wis{GL}_n$-orbits in $\gamma^{-1}(\tau)$. Such a point $\phi$ must be of the form
\[
x \mapsto 0, \quad y \mapsto 0, \quad z \mapsto 0, \quad X \mapsto K, \quad Y \mapsto L, \quad Z \mapsto M \]
By semi-stability, $(K,L,M)$ defines a simple $n$-dimensional representation of $A$ and its $\wis{GL}_n \times \C^*$-orbit defines the point $[det(K):det(L):det(M)] \in \mathbb{P}^2_c$. hence we may assume for instance that $K$ is invertible.

The tangent space $T_{\phi} \wis{rep}_n^{ss} B$ is the linear space of all trace-preserving algebra maps $B \rTo M_n(\C[\epsilon])$ of the form
\[
x \mapsto 0+\epsilon U, y \mapsto 0+\epsilon V, z \mapsto 0+\epsilon W, X \mapsto K + \epsilon R, Y \mapsto L + \epsilon S, Z \mapsto M + \epsilon T \]
and we have to use the relations in $B$ to show that the dimension of this space is at most $n^2+3$. As $(K,L,M)$ is a simple $n$-dimensional representation of the Sklyanin algebra, we know already that $(R,S,T)$ depend on at most $n^2+2$ parameters. Further, from the commutation relations in $B$ we deduce the following equalities (using the assumption that $K$ is invertible)
\begin{itemize}
\item{$xX = Xx \Rightarrow UK = KU$}
\item{$xY=Xy \Rightarrow UL=KV \Rightarrow K^{-1}UL = V$}
\item{$xZ=Xz \Rightarrow UM=KW \Rightarrow K^{-1}UM = W$}
\item{$Yx=yX \Rightarrow LU=VK \Rightarrow LK^{-1}U = V$}
\item{$Zx=zX \Rightarrow MU=WK \Rightarrow MK^{-1}U = W$}
\end{itemize}
These equalities imply that $K^{-1}U$ commutes with $K,L$ and $M$ and as $(K,L,M)$ is a simple representation and hence generate $M_n(\C)$ it follows that $K^{-1}U = \lambda 1_n$ for some $\lambda \in \C$. But then it follows that
\[
U = \lambda K, \quad V = \lambda L, \quad W = \lambda M \]
and so the triple $(U,V,W)$ depends on at most one extra parameter, showing that $T_{\phi} \wis{rep}^{ss}_n B$ has dimension at most $n^2+3$, finishing the proof. $\hfill \qed$

\begin{remark}
The statement of the previous theorem holds in a more general setting, that is, $\wis{rep}^{ss}_n B$ is smooth whenever $B = A \oplus A^+ t \oplus (A^+)^2t^2 \oplus \ldots$ with $A$ a positively graded algebra that is Azumaya away from the maximal ideal $A^+$ and $Z(A)$ smooth away from the origin. 
\end{remark}

\vskip 3mm

Unfortunately this does not imply that $\wis{proj} Z(B) = \wis{rep}^{ss}_n B // \wis{GL}_n \times \C^*$ is smooth as there are closed $\wis{GL}_n \times \C^*$ orbits with stabilizer subgroups strictly larger than $\C^* \times 1$. This happens precisely in semi-stable representations $\phi$ determined by
\[
x \mapsto 0, \quad y \mapsto 0, \quad z \mapsto 0, \quad X \mapsto K, \quad Y \mapsto L, \quad Z \mapsto M \]
with $[det(K):det(L):det(M)] \in E'$. In which case the matrices $(K,L,M)$ can be brought into the form
{\tiny
\[
\begin{bmatrix}
0 & 0 &  \hdots & \hdots & a_{n-1} t \\
a_0 & 0 &  \hdots & \hdots & 0 \\
0 & a_1 &  \ddots &  & \vdots \\
\vdots &  & \ddots & \ddots & \vdots \\
0 & 0  & \hdots & a_{n-2} & 0 
\end{bmatrix} \quad 
\begin{bmatrix}
0 & 0 &  \hdots & \hdots & b_{n-1} t \\
b_0 & 0 &  \hdots & \hdots & 0 \\
0 & b_1 &  \ddots &  & \vdots \\
\vdots &  & \ddots & \ddots & \vdots \\
0 & 0  & \hdots & b_{n-2} & 0 
\end{bmatrix}
\quad 
\begin{bmatrix}
0 & 0 &  \hdots & \hdots & c_{n-1} t \\
c_0 & 0 &  \hdots & \hdots & 0 \\
0 & c_1 &  \ddots &  & \vdots \\
\vdots &  & \ddots & \ddots & \vdots \\
0 & 0  & \hdots & c_{n-2} & 0 
\end{bmatrix}
\]
}
and the stabilizer subgroup is generated by $\C^* \times 1$ together with the cyclic group of order $n$
\[
\pmb{\mu}_n = \langle ( \begin{bmatrix} 1 & & & \\ & \zeta & & \\ & & \ddots & \\ & & & \zeta^{n-1} \end{bmatrix}, \zeta ) \rangle \]

\begin{lemma} If $\phi$ is a representation as above, then the normal space $N(\phi)$ to the $\wis{GL}_n$-orbit decomposes as a representation over the $\wis{GL}_n \times \C^*$-stabilizer subgroup $\pmb{\mu}_n$ as $\C_0 \oplus \C_0 \oplus \C_3 \oplus \C_{-1}$, that is, the associated local weighted quiver is 
\vskip 4mm
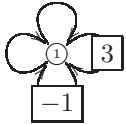
\[
\xymatrix{ \vtx{1} \ar@(dl,ul) \ar@(ul,ur) \ar@(ur,dr)|{\boxed 3} \ar@(dr,dl)|{\boxed {-1}}} \]
\end{lemma}

\begin{proof} The extra tangential coordinate $\lambda$ determines the tangent-vectors of the three degree zero generators
\[
x \mapsto 0 + \epsilon \lambda K, \qquad y \mapsto 0 + \epsilon \lambda L, \qquad z \mapsto 0 + \epsilon \lambda M \]
and so the generator of $\pmb{\mu}_n$ acts as follows
\[
\begin{bmatrix} 1 & & & \\ & \zeta^{n-1} & & \\ & & \ddots & \\ & & & \zeta \end{bmatrix}.(\epsilon \lambda (K,L,M)). \begin{bmatrix} 1 & & & \\ & \zeta & & \\ & & \ddots & \\ & & & \zeta^{n-1} \end{bmatrix} =\epsilon \zeta^{n-1} \lambda (K,L,M) \]
and hence accounts for the extra component $\C_{-1}$.
\end{proof}

We have now all information to prove Theorem 3 which asserts that the canonical map
\[
\wis{proj} Z(B) \rOnto \wis{spec} Z(A) \]
is a partial resolution of singularities, with singular locus $E' = E/\langle p \rangle$ in the exceptional fiber, all singularities of type $\C \times \C^2/\Z_n$. In other words, the isolated singularity of $\wis{spec} Z(A)$ `sees' the elliptic curve $E'$ and the isogeny $E \rOnto E'$ defining the $3$-dimensional Sklyanin algebra $A$.

\vskip 3mm
\noindent
{\bf Proof of Theorem 3 : } The GIT-quotient map
\[
\wis{rep}^{ss}_n B \rOnto \wis{proj} Z(B) \]
is a principal $\wis{PGL}_n \times \C^*$-bundle away from the elliptic curve $E'$ in the exceptional fiber whence $\wis{proj} Z(B) - E'$ is smooth. The application to the Luna slice theorem of \cite[Thm. 5]{BocklandtSymens} asserts that for any point $\overline{\phi} \in E' \rInto \wis{proj} Z(B)$ and all $t \in \C$ there is a neighborhood of $(\overline{\phi},t) \in \wis{proj} Z(B) \times \C$ which is \'etale isomorphic to a neighborhood of $0$ in $N(\phi) // \pmb{\mu}_n$. From the previous lemma we deduce that
\[
N(\phi) // \pmb{\mu}_n \simeq \C \times \C \times \C^2 / \mathbb{Z}_n \]
where $\C[\C^2/\mathbb{Z}_n] \simeq \C[u,v,w]/(w^n-uv^3)$, finishing the proof.

\vskip 3mm

As $B$ is a finite module over its center, it defines a coherent sheaf of algebras over $\wis{proj} Z(B)$. From Theorem 3 we obtain

\begin{corollary} The sheaf of Cayley-Hamilton algebras $\mathcal{B}$ on $\wis{proj} Z(B)$ is Azumaya away from the elliptic curve $E'$ in the exceptional fiber $\pi^{-1}(\mathfrak{m}) = \PP^2$ and hence is Cayley-smooth on this open set. However, $\mathcal{B}$ is not Cayley-smooth.
\end{corollary}
\begin{proof}
For a point $p$ in the exceptional fiber $\pi^{-1}(\mathfrak{m})-E'$ we already know that $\wis{proj} Z(B)$ is smooth and that $B$ is Azumaya, which implies that $\wis{rep}^{ss}_n B$ is smooth in the corresponding orbit. However, for a point $p \in E'$ we know that $\wis{proj} Z(B)$ has a non-isolated singularity in $p$. Therefore, $\wis{rep}_n^{ss}\mathcal{B}$ can not be smooth in the corresponding orbit, as the only central singularity possible for a Cayley-smooth order over a center of dimension 3 is the conifold singularity, which is isolated.
\end{proof}

\end{document}